\documentclass[english]{article}
\usepackage[T1]{fontenc}
\usepackage[latin9]{inputenc}
\usepackage{geometry}
\usepackage{amsmath}
\usepackage{amssymb}
\usepackage{setspace}
\makeatletter

\providecommand{\tabularnewline}{\\}

\newcommand{\lyxaddress}[1]{
\par {\raggedright #1
\vspace{1.4em}
\noindent\par}
}

\makeatother

\usepackage{babel}
\begin{document}

\title{Generalised Split Octonions and their transformation in SO(7) symmetry}

\author{K. Pushpa, P. S. Bisht and O. P. S. Negi}

\maketitle

\lyxaddress{\begin{center}
Department of Physics,\\
 Kumaun University, S. S. J. Campus, \\
Almora-263601 (Uttarakhand) India\\

\par\end{center}}

\lyxaddress{\begin{center}
Email-pushpakalauni60@yahoo.co.in\\
 ps\_bisht123@rediffmail.com \\
ops\_negi@yahoo.co.in.
\par\end{center}}
\begin{abstract}
Generators of $SO(8)$ group have been described by using direct product
of the Gamma matrices and the Pauli Sigma matrices. We have obtained
these generators in terms of generalized split octonion also. These
generators have been used to describe the rotational transformation
of vectors for $SO(7)$ symmetry group.
\end{abstract}

\section{Introduction}

Among all the division algebras, octonion forms the largest normed
division algebra \cite{key-1,key-2,key-3}. These are an algebraic
structure defined on the 8-dimensional real vector space and they
are non-associative and non-commutative extension of the algebra of
quaternions. G$\ddot{u}$naydin and G$\ddot{u}$rsey discussed the
Lie algebra of $G_{2}$ group and its embedding in $SO(7)$ and $SO(8)$
groups \cite{key-4}. A dynamical scheme of quark and lepton family
unification based on non associative algebra has also been discussed
\cite{key-5}. Generators of $SO(8)$ are constructed by using Octonion
structure tensors \cite{key-6} and the representations of these generators
are also given as the products of Octonions. Lassig and Joshi \cite{key-7}
introduced the bi-modular representation of octonions and formulated
the $SO(8)$ gauge theory equivalent to the octonionic construction,
also the peculiarities of the eight dimensional space has been described
by Gamba \cite{key-8}. Octonionic description has been used to describe
the various applications in quantum chromo dynamics \cite{key-9},
in the study of representations of Clifford algebras \cite{key-10},
non associative quantum mechanics \cite{key-11}, in the study of
symmetry breaking \cite{key-12}, in the study of flavor symmetry
\cite{key-13}, in proposals of unified field theory models \cite{key-14},
unitary symmetry \cite{key-15}, octonionic gravity \cite{key-16}.
In the present study, we have described the generalized split octonionic
description to represent the 8-dimensional algebra. We have obtained
the generators of $SO(8)$ group with the help of generalized split
octonions as well as direct product of Gamma and Pauli Sigma matrices.
Since $SO(7)\,$$\subset$~$SO(8)$, therefore the presented representation
also gives the generators of $SO(7)$ symmetry group. These generators
of $SO(7)$ group have been used to describe the rotational transformation
in 7 dimensional space.

\section{$SO(8)$ Symmetry group generators by using Direct Product of Gamma
and Sigma Matrices}

$SO(8)$ represents the special orthogonal group of eight-dimensional
rotations having 28 generators. By using direct product of Pauli matrices
and Gamma matrices, we have constructed the eight dimensional representation
of $SO(8)$ Symmetry.\\
Dirac Gamma matrices $(\gamma^{0},\gamma^{1},\gamma^{2},\gamma^{3})$
are a set of conventional matrices with specific anti-commutation
relations and they generate a matrix representation of the Clifford
algebra $C\ell_{1,3}(R$).
\begin{gather}
\gamma^{0}=\left[\begin{array}{cc}
\sigma^{0} & 0\\
0 & -\sigma^{0}
\end{array}\right]\,\,\,\,\gamma^{j}=\left[\begin{array}{cc}
0 & \sigma^{j}\\
-\sigma^{j} & 0
\end{array}\right]\,\,\forall j=1,2,3,\label{eq:1}
\end{gather}
where $\sigma^{0}$ is a 2 $\times$ 2 unit matrix and $\sigma^{j}$
are Pauli matrices defined as
\begin{align}
\sigma^{0}=\left[\begin{array}{cc}
1 & 0\\
0 & 1
\end{array}\right],\,\,\, & \sigma^{1}=\left[\begin{array}{cc}
0 & 1\\
1 & 0
\end{array}\right],\,\,\,\,\sigma^{2}=\left[\begin{array}{cc}
0 & -i\\
i & 0
\end{array}\right];\,\,\sigma^{3}=\left[\begin{array}{cc}
1 & 0\\
0 & -1
\end{array}\right].\label{eq:2}
\end{align}
The product of Gamma matrices are given by
\begin{align}
\gamma^{4}= & i\gamma^{0}\gamma^{1}\gamma^{2}\gamma^{3}=\left[\begin{array}{cccc}
0 & 0 & 1 & 0\\
0 & 0 & 0 & 1\\
1 & 0 & 0 & 0\\
0 & 1 & 0 & 0
\end{array}\right].\label{eq:3}
\end{align}
 Here we intend to describe the generators of $SO(8)$ symmetry group
by using the direct product of Gamma matrices with Pauli sigma matrices.
For this, we have chosen the following representations,
\begin{align}
\beta_{0} & =I_{4}\otimes\sigma^{0};\nonumber \\
\beta_{1} & =-i\,\gamma^{4}\otimes\sigma^{2}=[b_{16}-b_{25}+b_{38}-b_{47}+b_{52}-b_{61}+b_{74}-b_{83}];\nonumber \\
\beta_{2} & =\sigma^{1}\otimes\gamma^{3}=[b_{17}-b_{28}-b_{35}+b_{46}+b_{53}-b_{64}-b_{71}+b_{72}];\nonumber \\
\beta_{3} & =\sigma^{1}\otimes\gamma^{1}=[b_{18}+b_{27}-b_{36}-b_{45}+b_{54}+b_{63}-b_{72}-b_{81}];\nonumber \\
\beta_{4} & =\gamma^{4}\otimes\sigma^{0}=[b_{11}+b_{22}+b_{33}+b_{44}-b_{55}-b_{66}-b_{77}-b_{88}];\nonumber \\
\beta_{5} & =-i\,\gamma^{3}\otimes\sigma^{2}=[b_{16}-b_{25}-b_{37}+b_{46}+b_{52}-b_{61}-b_{74}+b_{83}];\nonumber \\
\beta_{6} & =-i\,\gamma^{2}\otimes\sigma^{0}=[-b_{17}-b_{28}+b_{35}+b_{46}+b_{53}+b_{64}-b_{71}-b_{82}];\nonumber \\
\beta_{7} & =-i\,\gamma^{^{1}}\otimes\sigma^{2}=[b_{18}-b_{27}+b_{36}-b_{45}-b_{54}+b_{63}-b_{72}+b_{81}];\label{eq:4}
\end{align}
where $I_{4}$ is 4$\times$4 identity matrix, $\beta_{0}$ is the
8$\times$8 identity matrix and $\beta_{a}$, (a=1,...7) represents
8 $\times$ 8 matrices. $b_{ij}$ represents the 8 $\times$ 8 matrices
in which $ij{}^{th}$ matrix element is unity and rest elements are
zero. $\beta_{a}$ represents the 7 generators of $SO(8)$ symmetry
group. It also satisfies the anti-commutation relation as
\begin{align}
\beta_{a}\beta_{b}+\beta_{b}\beta_{a}= & 2\delta_{ab}.\label{eq:5}
\end{align}
Furthermore, the rest 21 generators of $SO(8)$ symmetry group are
obtained by the following relation
\begin{align}
\beta_{ab}=\frac{1}{2i} & [\beta_{a},\beta_{b}].\label{eq:6}
\end{align}
Thus the 7 generators $\beta_{a}$ and 21 generators $\beta_{ab}$
form the total 28 generators of the $SO(8)$ symmetry group. In this
way, we have represented generators of $SO(8)$ symmetry group by
using the direct product of Gamma matrices and Pauli matrices. In
the succeeding sections, we show that these generators could also
been obtained with the help of generalized split octonions.

\section{Split Octonions}

Split Octonion algebra \cite{key-1} with its split base units is
defined as 
\begin{gather}
u_{0}=\frac{1}{2}\left(e_{0}+ie_{7}\right),\,\,\,\,\,\, u_{0}^{\star}=\frac{1}{2}\left(e_{0}-ie_{7}\right);\nonumber \\
u_{j}=\frac{1}{2}\left(e_{j}+ie_{j+3}\right)\,\,\,\,\, u_{m}^{\star}=\frac{1}{2}\left(e_{j}-ie_{j+3}\right).\label{eq:7}
\end{gather}
where j=1,2,3. \\
These basis elements satisfy the following algebra
\begin{align}
u_{i}u_{j}=-u_{j}u_{i}=\epsilon_{ijk}u_{k}^{\star},\,\,\,\,\,\, & u_{i}^{\star}u_{j}^{\star}=-u_{j}^{\star}u_{i}^{\star}=\epsilon_{ijk}u_{k};\nonumber \\
u_{i}u_{j}^{\star}=-\delta_{ij}u_{0},\,\,\,\,\,\,\,\,\,\,\,\,\,\,\,\, & u_{i}^{\star}u_{j}=-\delta_{ij}u_{0}^{\star};\nonumber \\
u_{0}u_{i}=u_{i}u_{0}^{\star}=u_{i},\,\,\,\,\,\,\,\,\,\,\,\,\, & u_{0}^{\star}u_{i}^{\star}=u_{i}^{\star}u_{0}=u_{i}^{\star};\nonumber \\
u_{i}u_{0}=u_{0}u_{i}^{\star}=0,\,\,\,\,\,\,\,\,\,\,\,\,\, & u_{i}^{\star}u_{0}^{\star}=u_{0}^{\star}u_{i}=0;\nonumber \\
u_{0}u_{0}^{\star}=u_{0}^{\star}u_{0}=0,\,\,\,\,\,\,\,\,\,\,\,\,\, & u_{0}^{2}=u_{0},\,\, u_{0}^{\star2}=u_{0}^{\star}.\label{eq:8}
\end{align}
These relations $\left(\ref{eq:8}\right)$ are invariant under $G_{2}$
transformation. Split octonion multiplication has been shown in Table
1.

\begin{table}[h]
\begin{centering}
\begin{tabular}{|c|c|c|c|c|c|c|c|c|}
\hline 
. & $u_{0}^{\star}$ & $u_{1}^{\star}$ & $u_{2}^{\star}$ & $u_{3}^{\star}$ & $u_{0}$ & $u_{1}$ & $u_{2}$ & $u_{3}$\tabularnewline
\hline 
$u_{0}^{\star}$ & $u_{0}^{\star}$ & $u_{1}^{\star}$ & $u_{2}^{\star}$ & $u_{3}^{\star}$ & 0 & 0 & 0 & 0\tabularnewline
\hline 
$u_{1}^{\star}$ & 0 & 0 & $u_{3}$ & -$u_{2}$ & $u_{1}^{\star}$ & $-u_{0}^{\star}$ & 0 & 0\tabularnewline
\hline 
$u_{2}^{\star}$ & 0 & $-u_{3}$ & 0 & $u_{1}$ & $u_{2}^{\star}$ & 0 & -$u_{0}^{\star}$ & 0\tabularnewline
\hline 
$u_{3}^{\star}$ & 0 & $u_{2}$ & -$u_{1}$ & 0 & $u_{3}^{\star}$ & $0$ & $0$ & $-u_{0}^{\star}$\tabularnewline
\hline 
$u_{0}$ & 0 & 0 & $0$ & $0$ & $u_{0}$ & $u_{1}$ & $u_{2}$ & $u_{3}$\tabularnewline
\hline 
$u_{1}$ & $u_{1}$ & -$u_{0}$ & 0 & 0 & 0 & 0 & $u_{3}^{\star}$ & $-u_{2}^{\star}$\tabularnewline
\hline 
$u_{2}$ & $u_{2}$ & $0$ & -$u_{0}$ & 0 & $0$ & -$u_{3}^{\star}$ & 0 & $u_{1}^{\star}$\tabularnewline
\hline 
$u_{3}$ & $u_{3}$ & $0$ & 0 & -$u_{0}$ & 0 & $u_{2}^{\star}$ & -$u_{1}^{\star}$ & $0$\tabularnewline
\hline 
\end{tabular}
\par\end{centering}

\caption{Split Octonion multiplication table}
\end{table}

\section{Generalized Split Octonions}

Let us consider a spinor $U$ generated from split octonions $u_{k}$
and conjugate of split octonions $u_{k}^{\star}$ as
\begin{align}
U= & \left(\begin{array}{c}
u_{k}\\
u_{k}^{\star}
\end{array}\right)\label{eq:9}
\end{align}
where $k=0,1,2,3.$ \\
Since split octonion represents the non associative and non division
algebra. It cannot be represented in terms of matrices. However, there
exists some real matrices \cite{key-2,key-4} which are associated
with split octonions that can be obtained as 
\begin{align}
u_{0}\left[\begin{array}{c}
\begin{array}{c}
u_{0}\\
u_{1}\\
u_{2}\\
u_{3}\\
u_{0}^{\star}\\
u_{1}^{\star}\\
u_{2}^{\star}\\
u_{3}^{\star}
\end{array}\end{array}\right]=\left[\begin{array}{c}
\begin{array}{c}
u_{0}\\
u_{1}\\
u_{2}\\
u_{3}\\
0\\
0\\
0\\
0
\end{array}\end{array}\right] & \left[\begin{array}{cccccccc}
1 & 0 & 0 & 0 & 0 & 0 & 0 & 0\\
0 & 1 & 0 & 0 & 0 & 0 & 0 & 0\\
0 & 0 & 1 & 0 & 0 & 0 & 0 & 0\\
0 & 0 & 0 & 1 & 0 & 0 & 0 & 0\\
0 & 0 & 0 & 0 & 0 & 0 & 0 & 0\\
0 & 0 & 0 & 0 & 0 & 0 & 0 & 0\\
0 & 0 & 0 & 0 & 0 & 0 & 0 & 0\\
0 & 0 & 0 & 0 & 0 & 0 & 0 & 0
\end{array}\right]\left[\begin{array}{c}
\begin{array}{c}
u_{0}\\
u_{1}\\
u_{2}\\
u_{3}\\
u_{0}^{\star}\\
u_{1}^{\star}\\
u_{2}^{\star}\\
u_{3}^{\star}
\end{array}\end{array}\right].\label{eq:10}
\end{align}
where 8$\times$8 matrix term in the right side represents the left
operator for split octonion basis element $u_{0}$ and we will denote
it as $U_{L0}$ in our further calculations and equation $(\ref{eq:10})$
could be written as
\begin{align}
u_{0}[U]= & U_{L0}[U].\label{eq:11}
\end{align}
Similarly we have calculated four 8$\times$8 matrices corresponding
to split octonions and four 8$\times$8 matrices corresponding to
conjugate of split octonions with the help of the split octonion multiplication
table 1.\\
Left multiplication of split octonion units by $U$ gives
\begin{align}
U_{L0}= & [a_{11}+a_{22}+a_{33}+a_{44}];\nonumber \\
U_{L1}= & [a_{38}-a_{47}+a_{52}-a_{61}];\nonumber \\
U_{L2}= & [-a_{28}+a_{46}+a_{53}-a_{71}];\nonumber \\
U_{L3}= & [a_{27}-a_{36}+a_{54}-a_{81}].\label{eq:12}
\end{align}
Left multiplication of conjugate of split octonions with $U$ are
given as
\begin{align}
U_{L0}^{\star}= & [a_{44}+a_{55}+a_{66}+a_{77}];\nonumber \\
U_{L1}^{\star}= & [a_{16}-a_{25}+a_{64}-a_{74}];\nonumber \\
U_{L2}^{\star}= & [a_{17}-a_{35}-a_{64}+a_{82}];\nonumber \\
U_{L3}^{\star}= & [a_{18}-a_{45}+a_{63}-a_{72}],\label{eq:13}
\end{align}
where $U_{L0}$, $U_{L1}$, $U_{L2}$, $U_{L3}$, $U_{L0}^{\star}$,
$U_{L1}^{\star}$, $U_{L2}^{\star}$, $U_{L3}^{\star}$ represent
the 8$\times$8 matrix representations named generalized split octonions
and $a_{ij}$ represents the 8$\times$8 matrix whose $ij^{th}$ element
is 1 and other terms are zero. By using the different combinations
of equation (\ref{eq:12}) and (\ref{eq:13}), we construct the following
8$\times$8 matrices

\begin{gather}
U_{0}=U_{L0}+U_{L0}^{\star}=\left[\begin{array}{cccc}
\sigma_{0} & 0 & 0 & 0\\
0 & \sigma_{0} & 0 & 0\\
0 & 0 & \sigma_{0} & 0\\
0 & 0 & 0 & \sigma_{0}
\end{array}\right];\,\,\,\,\, U_{1}=U_{L1}+U_{L1}^{\star}=i\left[\begin{array}{cccc}
0 & 0 & \sigma_{2} & 0\\
0 & 0 & 0 & \sigma_{2}\\
\sigma_{2} & 0 & 0 & 0\\
0 & \sigma_{2} & 0 & 0
\end{array}\right];\nonumber \\
U_{2}=U_{L2}+U_{L2}^{\star}=\left[\begin{array}{cccc}
0 & 0 & 0 & \sigma_{3}\\
0 & 0 & -\sigma_{3} & 0\\
0 & \sigma_{3} & 0 & 0\\
-\sigma_{3} & 0 & 0 & 0
\end{array}\right];\,\,\,\,\, U_{3}=U_{L3}+U_{L3}^{\star}=\left[\begin{array}{cccc}
0 & 0 & 0 & \sigma_{1}\\
0 & 0 & -\sigma_{1} & 0\\
0 & \sigma_{1} & 0 & 0\\
-\sigma_{1} & 0 & 0 & 0
\end{array}\right];\nonumber \\
U_{4}=U_{L0}-U_{L0}^{\star}=\left[\begin{array}{cccc}
\sigma_{0} & 0 & 0 & 0\\
0 & \sigma_{0} & 0 & 0\\
0 & 0 & -\sigma_{0} & 0\\
0 & 0 & 0 & -\sigma_{0}
\end{array}\right];\,\,\,\,\, U_{5}=U_{L1}-U_{L1}^{\star}=i\left[\begin{array}{cccc}
0 & 0 & \sigma_{2} & 0\\
0 & 0 & 0 & -\sigma_{2}\\
\sigma_{2} & 0 & 0 & 0\\
0 & -\sigma_{2} & 0 & 0
\end{array}\right];\nonumber \\
U_{6}=U_{L2}-U_{L2}^{\star}=\left[\begin{array}{cccc}
0 & 0 & 0 & -\sigma_{0}\\
0 & 0 & \sigma_{0} & 0\\
0 & \sigma_{0} & 0 & 0\\
-\sigma_{0} & 0 & 0 & 0
\end{array}\right];\,\,\,\,\, U_{7}=U_{L3}-U_{L3}^{\star}=i\left[\begin{array}{cccc}
0 & 0 & 0 & \sigma_{2}\\
0 & 0 & \sigma_{2} & 0\\
0 & -\sigma_{2} & 0 & 0\\
-\sigma_{2} & 0 & 0 & 0
\end{array}\right].\label{eq:14}
\end{gather}
Here $U_{A}$ (A=1,2,...7) gives the seven generators of $SO(8)$
symmetry group and the other 21 generators of $SO(8)$ group could
be obtained by taking the commutation relation of these matrices i.e.
\begin{align}
U_{AB}=\frac{1}{2i} & [U_{A},U_{B}]\,\,\forall A,B=1,2,....7.\label{eq:15}
\end{align}
$U_{A}$ and $U_{AB}$ give the 28 generators of $SO(8)$ group. These
28 generators of $SO(8)$ group are the same as those obtained by
the direct product of sigma and Gamma matrices. Now we would find
the connection between generalized split octonions and direct product
of sigma and Gamma matrices. From equation (\ref{eq:4}) and equation
(\ref{eq:14}), we have
\begin{align}
\beta_{1} & \Longrightarrow-i\,\gamma^{4}\otimes\sigma^{2}=U_{1}=U_{L1}+U_{L1}^{\star};\nonumber \\
\beta_{2} & \Longrightarrow\sigma^{1}\otimes\gamma^{3}=U_{2}=U_{L2}+U_{L2}^{\star};\nonumber \\
\beta_{3} & \Longrightarrow\sigma^{1}\otimes\gamma^{1}=U_{3}=U_{L3}+U_{L3}^{\star};\nonumber \\
\beta_{4} & \Longrightarrow\gamma^{4}\otimes\sigma^{0}=U_{4}=U_{L0}-U_{L0}^{\star};\nonumber \\
\beta_{5} & \Longrightarrow-i\,\gamma^{3}\otimes\sigma^{2}=U_{5}=U_{L1}-U_{L1}^{\star};\nonumber \\
\beta_{6} & \Longrightarrow-i\,\gamma^{2}\otimes\sigma^{0}=U_{6}=U_{L2}-U_{L2}^{\star};\nonumber \\
\beta_{7} & \Longrightarrow-i\,\gamma^{^{1}}\otimes\sigma^{2}=U_{7}=U_{L3}-U_{L3}^{\star}.\label{eq:16}
\end{align}
These 21 generators $U_{AB}$ represent the generators of $SO(7)$
group since $SO(7)\subset SO(8)$. Furthermore, with the help of these
21 generators of $SO(7)$ group, we would like to describe the rotational
transformation in seven dimensional space.

\section{Rotational transformation for SO(7) symmetry}

Let us consider a general spinor $\psi$ in seven dimensional space.
Under this symmetry, the spinor $\psi$ transforms as 
\begin{align}
\psi\longmapsto\psi^{\shortmid}= & \exp\left[\sum_{A=1}^{7}f_{A}U_{A}\right]\psi\nonumber \\
= & e^{X}\psi\label{eq:17}
\end{align}
where vector $X$ is, 
\begin{align}
X= & \sum_{A=1}^{7}f_{A}U_{A}.\label{eq:18}
\end{align}
with $f_{1}$, $f_{2}$, $f_{3}$,..........$f_{7}$ as the components
of the vector. On expanding $X$,
\begin{align}
X= & \left[\begin{array}{cccccccc}
f_{4} & 0 & 0 & 0 & 0 & if_{1}-f_{5} & if_{2}-f_{6} & if_{3}-f_{7}\\
0 & f_{4} & 0 & 0 & f_{5}-if_{1} & 0 & if_{3}+f_{7} & -if_{2}-f_{6}\\
0 & 0 & f_{4} & 0 & f_{6}-if_{2} & -if_{3}-f_{7} & 0 & if_{1}+f_{5}\\
0 & 0 & 0 & f_{4} & f_{7}-if_{3} & if_{2}+f_{6} & -if_{1}-f_{5} & 0\\
0 & if_{1}+f_{5} & if_{2}+f_{6} & if_{3}-f_{7} & -f_{4} & 0 & 0 & 0\\
-if_{1}-f_{5} & 0 & -if_{3}-f_{7} & f_{6}-if_{2} & 0 & -f_{4} & 0 & 0\\
if_{2}-f_{6} & f_{7}-if_{3} & 0 & if_{1}-f_{5} & 0 & 0 & -f_{4} & 0\\
if_{3}-f_{7} & if_{2}-f_{6} & f_{5}-if_{1} & 0 & 0 & 0 & 0 & -f_{4}
\end{array}\right]\label{eq:19}
\end{align}
which is a traceless Hermitian matrix. In compact form X can be written
as,
\begin{align}
X= & \left[\begin{array}{cc}
A & B^{\dagger}\\
B & -A
\end{array}\right]\label{eq:20}
\end{align}
where $A$ and $B$ are given as,
\begin{gather}
A=\left[\begin{array}{cccc}
f_{4} & 0 & 0 & 0\\
0 & f_{4} & 0 & 0\\
0 & 0 & f_{4} & 0\\
0 & 0 & 0 & f_{4}
\end{array}\right];\label{eq:21}\\
B=\left[\begin{array}{cccc}
0 & if_{1}-f_{5} & if_{2}-f_{6} & if_{3}-f_{7}\\
f_{5}-if_{1} & 0 & if_{3}+f_{7} & -if_{2}-f_{6}\\
f_{6}-if_{2} & -if_{3}-f_{7} & 0 & if_{1}+f_{5}\\
f_{7}-if_{3} & if_{2}+f_{6} & -if_{1}-f_{5} & 0
\end{array}\right].\label{eq:22}
\end{gather}
In equation (\ref{eq:20}), the term $B$ corresponds to the split
octonions ($u_{1}$,~$u_{2}$,~$u_{3}$) and $B^{\dagger}$ corresponds
to conjugate of split octonions ($u_{1}^{\star}$,~$u_{2}^{\star}$,$\, u_{3}^{\star}$).
Matrix $A$ represents the unit split octonions $u_{0}$ and $-A$
corresponds to $u_{0}^{\star}$. Matrices A and B are independent
of each other and they correspond to $u_{0}$ and $u_{j}$ respectively.\\
Furthermore, the constructed 8$\times$8 matrices given in equation
$\left(\ref{eq:16}\right)$ are being used to describe the rotation
in $SO(7)$ Symmetry. As an infinitesimal rotation by an angle $\theta$
in the plane ($k$, $l$) is obtained by the following operator \cite{key-1},
\begin{align}
R_{kl} & =1+\theta U_{k}U_{l}\label{eq:23}
\end{align}
which acts on a vector $X$ to form a rotated vector $X^{\shortmid}$
as,
\begin{align}
X^{\shortmid}= & R_{kl}XR_{kl}^{-1}\label{eq:24}
\end{align}
By using equation $\left(\ref{eq:23}\right)$, the rotation operator
$R_{12}$ becomes,
\begin{align}
R_{12}= & \left[\begin{array}{cccccccc}
1 & 0 & 0 & i\theta & 0 & 0 & 0 & 0\\
0 & 1 & i\theta & 0 & 0 & 0 & 0 & 0\\
0 & -i\theta & 1 & 0 & 0 & 0 & 0 & 0\\
-i\theta & 0 & 0 & 1 & 0 & 0 & 0 & 0\\
0 & 0 & 0 & 0 & 1 & 0 & 0 & i\theta\\
0 & 0 & 0 & 0 & 0 & 1 & i\theta & 0\\
0 & 0 & 0 & 0 & 0 & -i\theta & 1 & 0\\
0 & 0 & 0 & 0 & -i\theta & 0 & 0 & 1
\end{array}\right].\label{eq:25}
\end{align}
The rotation $R_{12}$ gives the rotated vector $X^{\shortmid}$ as
follows, 
\begin{align}
X^{\shortmid}= & X+2\theta\left[\begin{array}{cccccccc}
0 & 0 & 0 & 0 & 0 & if_{2} & -if_{1} & 0\\
0 & 0 & 0 & 0 & -if_{2} & 0 & 0 & if_{1}\\
0 & 0 & 0 & 0 & if_{1} & 0 & 0 & if_{2}\\
0 & 0 & 0 & 0 & 0 & -if_{1} & -if_{2} & 0\\
0 & if_{2} & -if_{1} & 0 & 0 & 0 & 0 & 0\\
-if_{2} & 0 & 0 & if_{1} & 0 & 0 & 0 & 0\\
if_{1} & 0 & 0 & if_{2} & 0 & 0 & 0 & 0\\
0 & -if_{1} & -if_{2} & 0 & 0 & 0 & 0 & 0
\end{array}\right].\label{eq:26}
\end{align}
While transforming these matrices in terms of rotations, we find that
the upper half of transformation corresponds to $u_{j}$ and lower
half corresponds to $u_{j}^{\star}$. So here we will be calculating
the rotational transformation of vectors by using split octonions
$u_{j}$ and Hermitian conjugate of split octonions $u_{j}^{\star}$
respectively. Since $\theta$ is an infinitesimal rotation, therefore
neglecting $\theta^{2}$ terms. Rotation $R_{12}$ transforms the
components of a vector $X$ as follows:
\begin{gather}
f_{1}\longmapsto f_{1}+2\theta[f_{2}-if_{5}];\nonumber \\
f_{2}\longmapsto f_{2}+2\theta[-f_{1}-if_{6}].\label{eq:27}
\end{gather}
We have calculated all the possible rotations of $R_{kl}$ and also
calculated the corresponding transformation of all combinations $kl$.
These are shown in table 2.

\begin{table}[h]
\begin{centering}
\begin{tabular}{|c|c|c|c|c|c|c|c|}
\hline 
 & $f_{1}$ & $f_{2}$ & $f_{3}$ & $f_{4}$ & $f_{5}$ & $f_{6}$ & $f_{7}$\tabularnewline
\hline 
\hline 
$R_{12}$ & $f_{2}-if_{5}$ & $-f_{1}-if_{6}$ &  &  &  &  & \tabularnewline
\hline 
$R_{13}$ & $f_{3}-if_{5}$ &  & $-f_{1}-if_{7}$ &  &  &  & \tabularnewline
\hline 
$R_{14}$ & $f_{4}-if_{5}$ &  &  & $-f_{1}$ &  &  & \tabularnewline
\hline 
$R_{15}$ & $-if_{1}$ &  &  &  & $-if_{5}$ &  & \tabularnewline
\hline 
$R_{16}$ & $f_{6}-if_{5}$ &  &  &  &  & $-f_{1}+if_{2}$ & \tabularnewline
\hline 
$R_{17}$ & $f_{7}-if_{5}$ &  &  &  &  &  & $-f_{1}+if_{3}$\tabularnewline
\hline 
$R_{23}$ &  & $f_{3}-if_{6}$ & $-f_{2}-if_{7}$ &  &  &  & \tabularnewline
\hline 
$R_{24}$ &  & $f_{4}-if_{6}$ &  & $-f_{2}$ &  &  & \tabularnewline
\hline 
$R_{25}$ &  & $f_{5}-if_{6}$ &  &  & $-f_{2}+if_{1}$ &  & \tabularnewline
\hline 
$R_{26}$ &  & $-if_{2}$ &  &  &  & $-if_{6}$ & \tabularnewline
\hline 
$R_{27}$ &  & $f_{7}-if_{6}$ &  &  &  &  & $-f_{2}+if_{3}$\tabularnewline
\hline 
$R_{34}$ &  &  & $f_{4}-if_{7}$ & $-f_{3}$ &  &  & \tabularnewline
\hline 
$R_{35}$ &  &  & $f_{5}-if_{7}$ &  & $-f_{3}+if_{1}$ &  & \tabularnewline
\hline 
$R_{36}$ &  &  & $f_{6}-if_{7}$ &  &  & $-f_{3}+if_{2}$ & \tabularnewline
\hline 
$R_{37}$ &  &  & $-if_{3}$ &  &  &  & $-if_{7}$\tabularnewline
\hline 
$R_{45}$ &  &  &  & $f_{5}$ & $-f_{4}+if_{1}$ &  & \tabularnewline
\hline 
$R_{46}$ &  &  &  & $f_{6}$ &  & $-f_{4}+if_{2}$ & \tabularnewline
\hline 
$R_{47}$ &  &  &  & $f_{7}$ &  &  & $-f_{4}+if_{3}$\tabularnewline
\hline 
$R_{56}$ &  &  &  &  & $f_{6}+if_{1}$ & $-f_{5}+if_{2}$ & \tabularnewline
\hline 
$R_{57}$ &  &  &  &  & $f_{7}+if_{1}$ &  & $-f_{5}+if_{3}$\tabularnewline
\hline 
$R_{67}$ &  &  &  &  &  & $f_{7}+if_{2}$ & $-f_{6}+if_{3}$\tabularnewline
\hline 
\end{tabular}
\par\end{centering}

\caption{Rotational transformation of SO(7) symmetry group}
\end{table}
The blank spaces in table 2 show that these components are unchanged
under the transformations. We have seen that rotational transformations
of the vector $f_{4}$ corresponding to $R_{14}$, $R_{24}$, $R_{34}$
rotations give $-f_{1}$, $-f_{2}$, $-f_{3}$ and the transformation
of the same corresponding to $R_{45}$, $R_{46}$, $R_{47}$ gives
$f_{5}$, $f_{6}$, $f_{7}$. This transformation represents rotations
corresponding to $SO(7)$ group corresponding to split octonions and
conjugate of split octonions. These transformation are obtained by
the combinations of $U_{A}$(A=1,2,3,5,6,7) with $U_{4}=U_{L0}-U_{L0}^{\star}$.
Therefore, these transformations are corresponding to $u_{0}$ and
$u_{0}^{\star}$ respectively. Rest of the transformation have been
obtained by different combinations of $U_{A}$ with $U_{B}$ (A, B
= 1,2,3,5,6,7) except $U_{4}$. These transformations are obtained
by using the split octonion $u_{j}$ and Hermitian conjugate of split
octonion $u_{j}^{\star}$ respectively.

\section{Discussion}

The generators of $SO(8)$ and $SO(7)$ groups have been generated
by using the split octonion as a spinor and the same relation has
been also described with the help of direct product of Gamma matrices
and Pauli sigma matrices. By using these generators we have obtained
the 21 rotational transformation in the $R_{kl}$ plane for $SO(7)$
group. An infinitesimal rotation transformations of $SO(7)$ group
is defined by using generators of $SO(8)$ symmetry. These representations
define the generators of 7 and 8 dimensional orthogonal groups in
terms of split octonionic descriptions of the Clifford groups. Rotational
transformations in seven dimensional space by using direct product
of gamma matrices with generators of $SU(2)$ group has been described.
We have used these generators to describe the rotational transformation
in seven dimensions for $SO(7)$ symmetry corresponding to split octonion
$u_{0}$, $u_{j}$ and its Hermitian conjugate $u_{0}^{\star}$, $u_{j}^{\star}$
respectively. Dirac Gamma matrices $(\gamma^{0},\gamma^{1},\gamma^{2},\gamma^{3})$
are matrices with satisfy specific anti commutation relations that
ensure they generate a matrix representation of the Clifford algebra
$C\ell1,3(R)$. It is also possible to define higher-dimensional gamma
matrices by using the direct product of Dirac Gamma matrices with
$SU(N)$ group generator. This will make higher dimensional Clifford
algebra. Clifford algebra can provide infinitesimal spatial rotations
and Lorentz boosts. This generalized split octonion algebra can be
further used to describe the higher dimensional algebra. Also one
can calculate the rotational transformations and boosts in higher
dimensional algebra by using direct product of Gamma matrices with
generators of $SU(N)$ group.

\end{document}